\documentclass[12pt]{article}
\usepackage{amsmath}
\usepackage{amssymb}
\usepackage{dsfont}
\usepackage{cite}
\newsymbol \blackbox 1004
\newcommand{\eh}{\hfill}\newlength{\sperr}

\def\nn{\nonumber}

\def\vk{\varkappa}

\def\s{\sigma}
\def\la{\lambda}

\def\vp{\varphi}
\def\vt{\vartheta}

\def\wt{\widetilde}

\def\BC{{\mathbb C}}
\def\BR{{\mathbb R}}
\def\BN{{\mathbb N}}

\def\clr{\mathcal{R}}

\newcommand{\I}{\mathrm{i}}

\newtheorem{Pa}{Paper}[section]
\newtheorem{Tm}[Pa]{{\bf Theorem}}

\newtheorem{Cy}[Pa]{{\bf Corollary}}
\newtheorem{Rk}[Pa]{{\bf Remark}}

\newtheorem{Dn}[Pa]{{\bf Definition}}

\newenvironment{dedication}
        {\vspace{1ex}\begin{quotation}\begin{center}\begin{em}}    
        {\par\end{em}\end{center}\end{quotation}}

\title{Discrete Dirac systems on the semiaxis: rational reflection coefficients and Weyl functions}

\author{B.~Fritzsche, B.~Kirstein, I.Ya.~Roitberg and A.L.~Sakhnovich}

\date{}
\begin{document}
\maketitle

\begin{dedication} 
\end{dedication}

\begin{abstract}   We consider the cases of the self-adjoint and skew-self-adjoint discrete Dirac systems,
obtain explicit expressions for reflection coefficients and show that rational reflection coefficients
and Weyl functions coincide. 
\end{abstract}

{MSC(2010): 39A10, 39A12, 47A40}

\vspace{0.2em}

{\bf Keywords:}    Discrete self-adjoint  Dirac system, discrete skew-self-adjoint Dirac system, Weyl function,
reflection coefficient, B\"acklund-Darboux transformation.

\section{Introduction}\label{Intro}
\setcounter{equation}{0}
Discrete self-adjoint and skew-self-adjoint Dirac systems play an essential role in the study
of Toeplitz matrices (and corresponding measures), of discrete integrable nonlinear equations
(including isotropic Heisenberg magnet model) and of spectral theory of difference equations
(see, e.g., \cite{DerSi, FKKS, KaS, ALSJFA, ALSverb} and references therein). Weyl-Titchmarsh theory
of discrete systems is actively studied (see, e.g., \cite{CGe, Jo, SZ, ZeC} and various references therein).
In particular, Weyl--Titchmarsh theory of  discrete self-adjoint and skew-self-adjoint Dirac systems
was studied in \cite{FKKS, FKRS14, FKRS-LAA, KaS, RoS, SaSaR} (see also references
therein). It is known that Weyl--Titchmarsh (or simply Weyl) functions of continuous
Dirac systems  on the semi-axis are closely
related to the scattering data. Some particular results for the self-adjoint systems are contained, for instance, in \cite{BeM, GKS6}  
and the general cases of continuous self-adjoint and skew-self-adjoint systems were treated in the recent paper \cite{ALS-Scat}.
The present article may be considered as the continuation of the paper \cite{ALS-Scat}, where the important discrete case
is dealt with.
We consider the cases of the self-adjoint and skew-self-adjoint discrete Dirac systems,
obtain explicit expressions for reflection coefficients and show that rational reflection coefficients
and Weyl functions coincide.

General-type 
discrete self-adjoint Dirac system has the form:
\begin{equation} \label{0.1}
y_{k+1}(z)=(I_m+ \I z
j C_k)
y_k(z) \quad \left( k \in \BN_0
\right),
\end{equation}
where $\BN_0$ stands for the set of non-negative integers,  $I_m$ is the $m \times m$ identity matrix, $``\I"$ is the imaginary unit
($\I^2=-1$) and the $m \times m$ matrices $C_k$ are positive and  $j$-unitary:
\begin{equation} \label{0.2}
C_k>0, \quad C_k j C_k=j, \quad  j: = \left[
\begin{array}{cc}
I_{m_1} & 0 \\ 0 & -I_{m_2}
\end{array}
\right] \quad (m_1+m_2=m; \, \, m_1, \, m_2 > 0).
\end{equation}
We introduce the Jost solution and reflection coefficient of the system \eqref{0.1}, \eqref{0.2} in a similar to the continuous case way.
Namely, the Jost solution $\{F_k(z) \}$ ($z\in \BR$) of the Dirac system \eqref{0.1}, \eqref{0.2} is defined via its asymptotics
\begin{align} & \label{R!}
F_k(z)=\big(I_m +\I z
 j
\big)^{k}\big(I_m+o(1)\big), \quad k \to \infty.
\end{align}
The reflection coefficient $\clr(z)$ is introduced via the blocks of $F_0(z)$:
\begin{align} & \label{R1!}
\clr(z)= \begin{bmatrix} I_{m_1} & 0 \end{bmatrix} F_0(z)\begin{bmatrix} 0 \\ I_{m_2}  \end{bmatrix}\left(\begin{bmatrix} 0 & I_{m_2}  \end{bmatrix}F_0(z)
\begin{bmatrix} 0 \\ I_{m_2}  \end{bmatrix}\right)^{-1}.
\end{align}

Discrete skew-self-adjoint Dirac system (SkDDS) is given (see \cite{FKRS-LAA, KaS}) by the formula:
\begin{equation} \label{d1}
y_{k+1}(z)=\left(I_m+ \frac{\I}{ z}
C_k\right)
y_k(z),  \quad C_k=U_k^*jU_k \quad \left( k \in \BN_0
\right),
\end{equation}
where the matrices $U_k$ are unitary and $j$ is defined in \eqref{0.2}.

Direct and inverse problems (in terms of Weyl functions) were solved for systems \eqref{0.1}, \eqref{0.2} in \cite{FKRS14, RoS}
and for systems \eqref{d1} in \cite{FKKS, FKRS-LAA}. In particular,  in the case of rational Weyl matrix functions, direct and inverse
problems were solved explicitly using our GBDT version \cite{ALS94, ALS-JMAA, SaSaR} of the B\"acklund-Darboux transformation.
For various versions of B\"acklund-Darboux transformations and related
commutation methods see, for instance, \cite{Ci, D, GeT, Gu, KoSaTe, Mar, MS} and references therein.

The results of the paper imply that the procedures to recover systems from the Weyl functions enable us to recover
systems from the reflection coefficients as well. In the next section, we give some preliminary definitions and results
in order to make the paper self-sufficient. Two subsections of Section \ref{Reflect} are dedicated to
the reflection coefficients in the self-adjoint and skew-self-adjoint cases.

 In the paper, $\BN$ denotes the set of natural numbers, $\BR$ denotes the real axis, $\BC$ stands for the complex plane, and
$\BC_+$ ($\BC_-$) stands for the open upper (lower) half-plane. The spectrum of a square matrix $A$ is denoted by $\s(A)$.
\section{Preliminaries }\label{Prel}
\setcounter{equation}{0}
\subsection{Self-adjoint case}
\paragraph{1.} 
Self-adjoint discrete Dirac system and stability of the explicit procedure to recover it from the Weyl function
was studied in our recent paper \cite{RoS}. We refer to \cite{RoS} for the preliminary definitions and results in this subsection.
The fundamental $m \times m$
solution $\{W_k \}$ of  \eqref{0.1}  is normalized by
\begin{align} \label{0.3}&
W_0(z)=I_m.
\end{align}
\begin{Dn} \label{defWeyl}  The Weyl function of the Dirac system \eqref{0.1} 
$($which is given on the semi-axis $0\leq k < \infty$ and satisfies \eqref{0.2}$)$ is an $m_1\times m_2$ matrix  function $\vp(z)$ in the
lower half-plane, such that the following inequalities hold $:$
\begin{align} \label{0.4}&
\sum_{k=0}^\infty q(z)^k 
\begin{bmatrix}
 \vp(z)^* &I_{m_2} 
\end{bmatrix}
W_k(z)^*C_k W_k (z)
\begin{bmatrix}
 \vp(z) \\ I_{m_2} 
\end{bmatrix}<\infty \quad (z\in \BC_-),
\\ & \label{0.5}
q(z):=(1+|z|^2)^{-1}.
\end{align}
\end{Dn}
(For the case $z \in \BC_+$, the definition of the Weyl function $\vp(z)$  of Dirac system \eqref{0.1}, \eqref{0.2} 
was given in  \cite{FKRS14}.)

\paragraph{2.}  In order to consider the case of rational Weyl functions, we introduce generalized B\"acklund-Darboux transformation (GBDT)
of the discrete self-adjoint  Dirac systems. Each GBDT of the initial discrete Dirac system is determined by a triple $\{A, S_0, \Pi_0\}$ of  parameter matrices.
Here, we take a trivial initial system and
choose $n\in \BN$, two $n
\times n$
parameter matrices $A$ ($\det A \not=0$) and
$S_0>0$, and an $n
\times m$ parameter matrix $\Pi_0$ such that
\begin{equation} \label{0.6}
A S_0-S_0 A^*=\I \Pi_0 j \Pi_0^*.
\end{equation}
Define  the sequences $\{\Pi_r\}$ and $\{S_r\}$ $(r \geq 0)$  using the triple $\{A, S_0, \Pi_0\}$
and recursive
relations
\begin{align} & \label{0.7}
\Pi_{k+1}=\Pi_k+\I A^{-1}\Pi_k j \quad (k\geq 0),
\\ &  \label{0.8}
S_{k+1}=S_k+A^{-1}S_k (A^*)^{-1}+A^{-1}\Pi_k
\Pi_k^*(A^*)^{-1} \quad (k\geq 0).
\end{align}
From \eqref{0.6}--\eqref{0.8},  the validity of the matrix identity
\begin{equation} \label{0.9}
A S_{r}-S_{r} A^*=\I \Pi_{r} j \Pi_{r}^* \quad
(r \geq 0)
\end{equation}
follows by induction.  In the self-adjoint case, we introduce {\it admissible} triples $\{A, S_0, \Pi_0\}$ in the
following way.
\begin{Dn} \label{adm} The triple $\{A, S_0, \Pi_0\}$, where $\det A \not=0$, $S_0>0$ and 
\eqref{0.6} holds, is called admissible.
\end{Dn}
In view of \eqref{0.8}, for the admissible triple we have $S_r>0$ $(r\geq 0)$. Thus, the
sequence (potential) $\{C_k\}$ $(k\geq 0)$ is well-defined by the equality
\begin{equation} \label{0.10}
C_k:=I_m+\Pi_k^*S_k^{-1}\Pi_k-\Pi_{k+1}^*S_{k+1}^{-1}\Pi_{k+1} .
\end{equation}
Moreover, from \cite[Theorem 2.5]{RoS} we see that the matrices $C_k$ satisfy \eqref{0.2}.
We say that the potential $\{C_k\}$ is {\it determined} by the
admissible triple. The potential determined by an admissible triple $\{A, S_0, \Pi_0\}$ is called {\it pseudo-exponential}.
We note that the notion of the pseudo-exponential (and {\it strictly pseudo-exponential}) potentials for the self-adjoint continuous
case was introduced first in \cite{GKS1} (see also \cite{GKS6}). In the discrete case, some additional requirements on the
{\it admissible and strongly admissible triples} (which determine pseudo-exponential and strictly pseudo-exponential,
respectively, potentials) appear.

All Weyl functions $\vp(z)$ are contractive in $\BC_-$
and all the potentials $\{C_k\}$, such that $\vp(z)$ (for the corresponding systems) are contractive and $\vp(-1/z)$ are
strictly proper rational, are determined by some admissible triples \cite{RoS}. Strongly admissible triples for the
self-adjoint case are considered in Subsection \ref{Refl1}.

We will need also the matrix function $w_A$,
which for each $k \geq 0$  is
a so called
transfer matrix function in Lev Sakhnovich form
\cite{SaL1, SaL3, SaSaR} and is defined by the relation
\begin{equation} \label{0.11}
w_A(k,\lambda):=I_m-\I j \Pi_k^*S_k^{-1}(A-\lambda
I_n)^{-1}\Pi_k.
\end{equation}
The fundamental
solution
$\{W_{k}\}$ of the Dirac system \eqref{0.1}
admits the representation
\begin{equation} \label{0.12}
W_{k}(z)=w_A(k,\,-1/z)\big(I_m +\I z
 j
\big)^{k}w_A(0,\,-1/z)^{-1} \quad (k \geq 0),
\end{equation}
where $w_A$ is defined in \eqref{0.11}.

Now, we partition $\Pi_k$ and write it down in the form
\begin{align}  \label{v7}  &
\Pi_k=\begin{bmatrix}(I_n+\I A^{-1})^k\vt_1 & (I_n-\I A^{-1})^k \vt_2 \end{bmatrix},
\end{align} 
where $\vt_1$ and $\vt_2$ are $n\times m_1$ and $n\times m_2$, respectively, blocks of $\Pi_0$. Assume further in this subsection that
\begin{align}  \label{R0}  &
\pm \I \not\in \s(A).
\end{align}
In view of \eqref{0.9} and \eqref{v7}, setting
\begin{align}  \label{v8}  &
R_{r}:=(I_n+\I A^{-1})^{-r}S_r\big(I_n-\I (A^*)^{-1}\big)^{-r}
\end{align} 
we have
\begin{align} \nn
R_{k+1}-R_k=& 2(I_n+\I A^{-1})^{-k-1}A^{-1}(I_n-\I A^{-1})^k \vt_2\vt_2^* \big((I_n-\I A^{-1})^k\big)^*\big(A^{-1}\big)^*
\\  \label{v9}  &
\times
\big((I_n+\I A^{-1})^{-k-1}\big)^*\geq 0.
\end{align} 
Since $R_0=S_0>0$, relations \eqref{v9} imply that there is a
limit
\begin{align}  \label{v10}  &
\lim_{k\to \infty}R_k^{-1}=\vk_R\geq 0.
\end{align} 
In a similar way we introduce the matrices
\begin{align}  \label{v8'}  &
Q_{r}:=(I_n-\I A^{-1})^{-r}S_r\big(I_n+\I (A^*)^{-1}\big)^{-r} ,
\end{align} 
and show that
\begin{align} \label{v9'}  &
Q_{k+1}-Q_k \geq 0.
\end{align} 
Since $Q_0=S_0>0$, relations \eqref{v9'} imply that there is a
limit
\begin{align}  \label{v10'}  &
\lim_{k\to \infty}Q_k^{-1}=\vk_Q\geq 0.
\end{align} 
\subsection{Skew-self-adjoint case}
The preliminary definitions and results on the skew-self-adjoint discrete  Dirac systems (SkDDS) \eqref{d1}
we take from \cite{FKKS} and sometimes from \cite{FKRS-LAA}. 
\begin{Rk}\label{RkNn}
The notations here slightly differ
from the notations in \cite{FKKS, FKRS-LAA}. In particular, we introduce the matrices $R_k$ and $Q_k$
in the both self-adjoint and skew-self-adjoint cases via formulas \eqref{v8} and \eqref{v8'}, respectively,
but in \cite{FKRS-LAA} $R_k$ stands for $Q_k$ in the the present notations and $Q_k$ stands for $R_k$.
\end{Rk}
\begin{Dn}\label{DnD0} The Weyl function of SkDDS is an $m_1 \times m_2$ matrix function $\vp(z)$ in 
$$\BC_M=\{z\in \BC: \, \Im(z)>M\}  \quad {\mathrm{for \,\, some}} \quad  M>0,$$
which satisfies the inequality
\begin{equation} \label{k1}
\sum_{k=0}^{\infty}\begin{bmatrix}  \vp(z)^* & I_{m_2}\end{bmatrix}w_k(z)^* w_k(z)
\left[
\begin{array}{c}
 \vp(z) \\ I_{m_2}
\end{array}
\right] < \infty,
\end{equation}
where $w_k(z)$ is the {fundamental solution} of SkDDS normalized by $w_0(z) \equiv I_m$.
\end{Dn}
Let us fix again an integer $n>0$, and consider an $n \times n$ matrix $A$
with $\det \, A \not= 0$, an  $n \times n$ matrix $S_0>0$  and  an $n \times m$  matrix $\Pi_0$. These
matrices should satisfy the  identity
\begin{equation} \label{k2}
A S_0-S_0 A^*= \I \Pi_0 \Pi_0^*.
\end{equation}
The sequences $\{\Pi_k\}$, $\{S_k\}$ and $\{C_k\}$ $(k \geq 0)$ are introduced using the triple $\{A, \, S_0, \, \Pi_0\}$ and relations
\begin{align}\label{k3}
&\Pi_{k+1}= \Pi_k+\I A^{-1} \Pi_k j, \\
&
S_{k+1}=S_k+ A^{-1} S_k (A^*)^{-1}+ A^{-1} \Pi_k j \Pi_k^*
(A^*)^{-1},\label{k4}
\\ \label{k5} &
C_k=j+ \Pi_k^* S_k^{-1} \Pi_k - \Pi_{k+1}^* S_{k+1}^{-1} \Pi_{k+1}.
\end{align}
Similar to the self-adjoint case we write down $\Pi_k$ in the form
$$\Pi_k=\begin{bmatrix}(I_n+\I A^{-1})^k\vt_1 & (I_n-\I A^{-1})^k \vt_2 \end{bmatrix}.$$
If $S_0>0$,  the identity \eqref{k2} holds and the  pair $\{A, \, \vt_1\}$
is controllable, then according to \cite[Lemma 3.2]{FKKS} and   \cite[Proposition 3.6]{FKKS} we have
$\det \, A \not= 0$, $S_k>0$ and the matrices $C_k$ admit representation $C_k=U_k^*jU_k$ 
from \eqref{d1}. That is, the sequence $\{C_k\}$ is well-defined and the corresponding
system is a skew-self-adjoint Dirac system. In the skew-self-adjoint case, the triple $\{A, \, S_0, \, \Pi_0\}$,
such that $S_0>0$,  the identity \eqref{k2} holds and the  pair $\{A, \, \vt_1\}$
is controllable, is called {\it admissible}. The potential determined by this triple
is called {\it pseudo-exponential}.

Moreover, if $\vp(z)$ is a strictly proper rational $m_1\times m_2$ matrix function then
it is the Weyl function of some skew-self-adjoint Dirac system with the pseudo-exponential potential  
(see \cite[Theorem 4.2]{FKRS-LAA}). We will require additionally that $\I \not\in \s(A)$.
\begin{Dn} In the skew-self-adjoint case, the triple $\{A, \, S_0, \, \Pi_0\}$, where  $S_0>0$, the identity \eqref{k2} is valid, the pair $\{A, \, \vt_1\}$
is controllable and $\I \not\in \s(A)$, is called strongly admissible. The potentials determined by the strongly admissible triples are called strictly pseudo-exponential.
\end{Dn}
Note that
\cite[Proposition 4.8]{FKRS-LAA} implies that if $S_0>0$,  \eqref{k2} holds and $0,\, \I \not\in \s(A)$ then $S_k>0$, the matrices $C_k$
are well-defined and there is a strongly admissible triple which determines the same potential $\{C_k\}$ as $\{A, \, S_0, \, \Pi_0\}$.
The fundamental solution $w_k$ of SkDDS determined by the  strongly admissible triple $\{A, \, S_0, \, \Pi_0\}$ has the form
\begin{align}& \label{k6}
w_k( z )=w_{ A }(k, -z ) \left(
I_{m}+\frac{\I}{z}j \right)^{k} w_{ A }(0, - z
)^{-1},
\end{align}
whereas $w_A$ in the skew-self-adjoint case is given by
\begin{align}
 & \label{k7}
w_{ A }(k, \la ):=I_{m}- \I \Pi_k^{*}
S_k^{-1} ( A- \lambda I_{n} )^{-1} \Pi_k.
\end{align}
Taking into account Remark \ref{RkNn}, we see that \cite[Proposition 4.10]{FKRS-LAA} and   \cite[(4.34)]{FKRS-LAA} imply that
\begin{align}  \label{k8}  &
\lim_{k\to \infty}Q_k^{-1}=0; \quad \lim_{k\to \infty}\big(Q_k^{-1}\wt G(A)^k\vt_1\big)=0, \quad \wt G(A):=(A-\I I_n)^{-1}(A+\I I_n)
\end{align}
in the case of a strongly admissible triple $\{A, \, S_0, \, \Pi_0\}$.
\section{Reflection coefficients}\label{Reflect}
\setcounter{equation}{0}
\subsection{Reflection coefficients: self-adjoint case }\label{Refl1}
In this subsection, we express (via the triple   $\{A, \, S_0, \, \Pi_0\}$)
the Jost solution and reflection coefficient, which are the analogs of the corresponding functions in the continuous case.

Uniqueness of the solution of the inverse problem to recover system from the Weyl function (see \cite[Theorem 2.3]{RoS}) together
with Theorems 2.6 and 2.8 and Proposition 2.7 (all from \cite{RoS}) imply that without loss of generality
one can require that $\s(A )\subset (\BC_+\cup \BR)$. Further we use a stronger requirement
\begin{align}  \label{R1}  &
\s(A )\subset \BC_+, \quad  \I \not\in \s(A).
\end{align}
Following \cite{FKKS, GKS1}, we call the admissible triple satisfying \eqref{R1} {\it strongly admissible}
and  we introduce the class of the strictly pseudo-exponential potentials
$\{C_k\}$.
\begin{Dn}\label{sps} The potentials $\{C_k\}$ of the Dirac systems \eqref{0.1}, \eqref{0.2},
which are determined by the strongly admissible triples, are called strictly pseudo-exponential.
\end{Dn}
In view of \eqref{0.11}, \eqref{v7}, \eqref{v8} and \eqref{v8'}, we have a representation
\begin{align}\label{R2}
& w_{A}(k,-1/ z )
\\ \nn &=
I_m-\I zj
\begin{bmatrix}
\vt_1^*R_k^{-1}(I_n+zA)^{-1}\vt_1 &
\vt_1^*R_k^{-1}(I_n+zA)^{-1} G(A)^k\vt_2 
\\
\vt_2^*(G(A)^k)^*R_k^{-1} (I_n+zA)^{-1}\vt_1 &
\vt_2^*Q_k^{-1}(I_n+zA )^{-1}\vt_2
\end{bmatrix},
\end{align}
where
\begin{align}  \label{R3}  &
G(A)=(I_n+\I A^{-1})^{-1}(I_n-\I A^{-1}).
\end{align} 
Relations \eqref{R1} and \eqref{R3} yield
\begin{align}  \label{R4}  &
\s\big(G(A)\big) \subset \{\la: \, |\la |<1 \}.
\end{align} 
Hence, from \eqref{v10}, \eqref{v10'} and \eqref{R2} we derive
\begin{align}\label{R5}
& \lim_{k\to \infty} w_{A}(k,-1/ z )
=
\begin{bmatrix}
\chi_1(z) & 0
\\
0 & \chi_2(z)
\end{bmatrix}, 
\\ \label{R6} &
 \chi_1(z):=I_{m_1}-\I z \vt_1^*\vk_R (I_n+zA)^{-1}\vt_1, \,\, \chi_2(z):=I_{m_2}+\I z \vt_2^*\vk_Q (I_n+zA)^{-1}\vt_2.
\end{align}
According to \eqref{R!}, \eqref{0.12} and \eqref{R5}, the Jost solution $\{F_k\}$ is given by the equalities
\begin{align}  \label{R7}  &
F_k(z)=W_k(z)w_{A}(0,-1/ z )\begin{bmatrix}
\chi_1(z)^{-1} & 0
\\
0 & \chi_2(z)^{-1}
\end{bmatrix}.
\end{align} 
Partition $w_A$ into the blocks corresponding to the partitioning of $j$:
\begin{align}\label{R8}&
 w_A(0,\lambda)=\left[
\begin{array}{lr}
a(\lambda) & b(\lambda) \\ c(\lambda) & d(\lambda)
\end{array}
\right].
\end{align}
It was shown in the proof of \cite[Theorem 2.6]{RoS} (see \cite[(2.29)]{RoS}) that the Weyl function $\vp(z)$
of the system \eqref{0.1}, \eqref{0.2} (in $\BC_-$)
is given by the formula
\begin{align}\label{R9}&
\vp(z)=b(-1/z)d(-1/z)^{-1}.
 \end{align}
Relations  \eqref{R1!}, \eqref{0.3} and \eqref{R7}--\eqref{R9} imply the following theorem.
\begin{Tm} \label{Tm1} Let Dirac system \eqref{0.1}, \eqref{0.2}
be a system with the
strictly pseudo-exponential potential $\{C_k\}$. Then the 
Weyl function $\vp(z)$ is the unique analytic continuation of the
reflection coefficient $\clr(z)$ of this system. That is, the reflection
coefficient and the Weyl function are given by the same rational matrix function.
\end{Tm}
From \cite[Theorem 2.6]{RoS} and Theorem \ref{Tm1} we derive the following corollary.
\begin{Cy}\label{Cy2} Let the potential $\{C_k\}$ be determined by a strongly admissible triple $\{A, \, S_0, \, \Pi_0\}$.
Then, the reflection coefficient of the Dirac system \eqref{0.1}, \eqref{0.2} is given by the formula
\begin{equation} \label{R10}
\clr(z)=-\I z\vt_1^*S_0^{-1}(I_n+zA^{\times})^{-1}\vt_2,
\quad A^{\times}=A+\I \vt_2 \vt_2^*S_0^{-1}.
\end{equation}
\end{Cy}
\subsection{Reflection coefficients: skew-self-adjoint case }\label{Refl2}
In the skew-self-adjoint case, we define the reflection coefficient $\clr(z)$ in a slightly more general way than in
the self-adjoint case.  That is, we consider the matrix valued  $m_2 \times m$ solution $Y$ of the system \eqref{d1}:
\begin{equation} \label{R11}
Y_k(z)=\big(1-\frac{\I}{z}\big)^k\left(\begin{bmatrix} 0 \\ I_{m_2}\end{bmatrix}+o(1)\right), \quad k\to \infty ,
\end{equation}
and set
\begin{align} & \label{R12}
\clr(z)= \begin{bmatrix} I_{m_1} & 0 \end{bmatrix} Y_0(z)\left(\begin{bmatrix} 0 & I_{m_2}  \end{bmatrix}Y_0(z)\right)^{-1}.
\end{align}

In order to express $\clr(z)$ via a strongly admissible  triple $\{A, \, S_0, \, \Pi_0\}$, we derive from \eqref{v7}, \eqref{v8}, \eqref{v8'} and
\eqref{k7} the representation
\begin{align}& \label{R13}
 w_{A}(k,-z )\begin{bmatrix} 0 \\ I_{m_2}\end{bmatrix}=
\begin{bmatrix} 0 \\ I_{m_2}\end{bmatrix}
-\I 
\begin{bmatrix}
\vt_1^*\big(\wt G(A)^k\big)^*Q_k^{-1}(z I_n+A)^{-1}\vt_2 
\\
\vt_2^*Q_k^{-1}(z I_n+A )^{-1}\vt_2
\end{bmatrix},
\end{align}
where $\wt G$ is introduced in \eqref{k8}. Formulas \eqref{k8} and \eqref{R13} imply that
\begin{align}& \label{R14}
\lim_{k \to \infty} w_{A}(k,-z )\begin{bmatrix} 0 \\ I_{m_2}\end{bmatrix}=\begin{bmatrix} 0 \\ I_{m_2}\end{bmatrix}.
\end{align}
It follows from \eqref{k6}, \eqref{R11} and \eqref{R14} that
\begin{align}& \label{R15}
 Y_k(z)=\big(1-\frac{\I}{z}\big)^k w_{A}(k,-z )\begin{bmatrix} 0 \\ I_{m_2}\end{bmatrix}.
 \end{align}
Hence, after we take into account \eqref{R12} and (similar to the self-adjoint case) partition $w_A$ (as in  \eqref{R8}), we obtain
\begin{align}& \label{R16}
 \clr(z)=b(-z)d(-z)^{-1}.
 \end{align}
On the other hand, according to \cite[(3.24)]{FKKS} the Weyl function $\vp(z)$ of the system \eqref{d1} 
is also given by the right-hand side of \eqref{R16}.
Thus, the following theorem is proved.
\begin{Tm} \label{Tm2} Let Dirac system \eqref{d1}
be a system with the
strictly pseudo-exponential potential $\{C_k\}$. Then the 
Weyl function $\vp(z)$ is the analytic continuation of the
reflection coefficient $\clr(z)$ of this system. More precisely, the reflection
coefficient and the Weyl function are given by the same rational matrix function.
\end{Tm}
The next corollary follows from \cite[Theorem 3.8]{FKKS} and Theorem \ref{Tm2}.
\begin{Cy}\label{Cy5} Let the potential $\{C_k\}$ be determined by a strongly admissible triple $\{A, \, S_0, \, \Pi_0\}$.
Then, the reflection coefficient of the skew-self-adjoint Dirac system \eqref{d1} is given by the formula
\begin{equation} \label{R17}
\clr(z)=-\I \vt_1^*S_0^{-1}(zI_n+zA^{\times})^{-1}\vt_2,
\quad A^{\times}=A-\I \vt_2 \vt_2^*S_0^{-1}.
\end{equation}
\end{Cy}

\bigskip

\noindent{\bf Acknowledgments.}
 {The research   of  A.L. Sakhnovich was supported by the
Austrian Science Fund (FWF) under Grant  No. P29177.}


\newpage 

\begin{flushright}
B. Fritzsche,\\
Fakult\"at f\"ur Mathematik und Informatik,  Universit\"at Leipzig, \\
 Augustusplatz 10,  D-04009 Leipzig, Germany, \\
 e-mail: {\tt Bernd.Fritzsche@math.uni-leipzig.de}
 
\vspace{0.5em}

B. Kirstein, \\ 
Fakult\"at f\"ur Mathematik und Informatik,  Universit\"at Leipzig, \\
 Augustusplatz 10,  D-04009 Leipzig, Germany, \\
 e-mail: {\tt Bernd.Kirstein@math.uni-leipzig.de}
 
\vspace{0.5em}

I.Ya. Roitberg, \\
 e-mail: {\tt  	innaroitberg@gmail.com}

\vspace{0.5em} 

A.L. Sakhnovich,\\
Faculty of Mathematics,
University
of
Vienna, \\
Oskar-Morgenstern-Platz 1, A-1090 Vienna,
Austria, \\
e-mail: {\tt oleksandr.sakhnovych@univie.ac.at}

\end{flushright}

\end{document}